\magnification=\magstep1

\input amstex
\documentstyle{amsppt}
\hsize = 6.5 truein
\vsize = 9 truein
\topmatter

\title {Variations on a result of Granville about polynomial invariance} \endtitle

\author Labib Haddad and Charles Helou  \endauthor

\address {120 Rue de Charonne, 75011 Paris, France; email: labib.haddad\@wanadoo.fr} \endaddress
\address {Penn State Univ., 25 Yearsley Mill Rd, Media, PA 19063, USA; e-mail: cxh22\@psu.edu} \endaddress

\endtopmatter

\document

\NoBlackBoxes

\head{Abstract}\endhead { In 2004, Andrew Granville proved the following: A polynomial in variables $x,y,z,$  invariant under the map $\varphi : z \to -(x+y+z)$ is a polynomial in $b = z(x+y+z)$. We tell,  curtly, its history  and then give a variation on this theme.}

\bigpagebreak

\head{History} \endhead

\medpagebreak

In March 2004, the second author received an inquiry from a colleague about the factorization of some symmetric polynomials, with power monomials, in two or three variables. There was a statement about a factorization in $\Bbb Z[x, y]$ of the polynomial $F_p(x,y)= (x+y)^p - x^p -y^p$, for a prime number $p>3$, namely 
$$F_p(x,y) = pxy(x+y)\left( x^2+xy+y^2 \right)^e C_p(x,y),$$
where $e=1 \  \text {or}\ 2$ according to whether $p\equiv -1 \ \text {or}  \ 1 \pmod 6$, respectively, and $C_p(x,y) \in \Bbb Z[x, y]$ is a so-called Cauchy polynomial, described as a well-known elementary result. In fact, $F_p(x,y)$ is a homogeneous polynomial of degree $p$, which, upon division by $y^p$, gives a polynomial $K_p(t)$, of degree $p-1$ in $t=\frac xy$, namely 
$$\frac {F_p(x,y)}{y^p} = \left( \frac xy +1 \right)^p -  \left( \frac xy \right)^p -1 = (t+1)^p -t^p -1 = \sum_{k=1}^{p-1} \binom pk t^k = K_p(t) .$$
Clearly, $K_p(0)=K_p(-1)=0$ and $K'_p(0)=p, \ K'_p(-1)=-p \neq 0$, so $0$ and $1$ are simple roots of $K_p$, so that $K_p$ is divisible by $t(t+1)$ in $\Bbb Z[t]$, and also by $p$ since the binomial coefficients $\binom pk$ are divisible by $p$ for $1\leq k \leq p-1$. Moreover, letting $\zeta_3 = \exp \left( \frac {2 i\pi}3 \right)$ be a primitive cubic root of unity, we have 
$$K_p(\zeta_3)= (\zeta_3 +1)^p - \zeta_3^p -1 = -(1+\zeta_3^p+\zeta_3^{2p})=0$$
and 
$$K'_p(\zeta_3)= p((\zeta_3+1)^{p-1} -\zeta_3^{p-1}) = p\zeta_3^{p-1}(\zeta_3^{p-1} -1),$$
so that $K'_p(\zeta_3)=0$ if and only if $p\equiv 1 \pmod 3$, i.e., $p\equiv 1 \pmod 6$. Thus $\zeta_3$ is a root of $K_p(t)$, and it's a double root if and only if $p\equiv 1 \pmod 6$, i.e., $K_p(t)$ is divisible by $t^2+t+1$ if $p\equiv -1 \pmod 6$ and by $(t^2+t+1)^2$ if $p\equiv 1 \pmod 6$, in $\Bbb Z[t]$. Thus
$$K_p(t) = pt(t+1)(t^2+t+1)^e C_p(t),$$
where $$e= \cases 1, &\text {if $p\equiv -1 \pmod 6$} \\  2, &\text {if $p\equiv 1 \pmod 6$}, \endcases$$ 
and $C_p(t) \in \Bbb Z[t]$, of degree $p-3-2e$, is called a Cauchy polynomial, because it was introduced by A. L. Cauchy in 1839 and reconsidered in 1841 (\cite{2, 3, 5}). The Cauchy polynomials were studied by D. Mirimanoff in 1903 (\cite{4}), who conjectured their irreducibility over $\Bbb Q$, a still open conjecture. For more details, one can consult the 1989 paper of G. Terjanian (\cite{6}). This factorization of $K_p(t)$ implies the factorization of $F_p(x,y)$ cited above. To be accurate, $C_p(t)$ should be $C_p(t,1)$ in the factorization of $K_p(t)$, but by ``abuse of notation", we identified the two.
 
This leads to introduce the polynomial 
$$E_2 = \frac {(x+y)^p -x^p -y^p}{pxy(x+y)} = \left( x^2+xy+y^2 \right)^e C_p(x,y),$$ 
which lies in $Z[x,y]$, and, by analogy, in the case of 3 variables $x, y, z$, the polynomial 
$$E_3 = \frac {(x+y+z)^p -x^p -y^p -z^p}{p(x+y)(y+z)(z+x)},$$ 
which lies in $Z[x,y,z]$. Indeed, letting 
$$H_p(x,y,z)= (x+y+z)^p -x^p -y^p -z^p = \sum \Sb 0\leq j, k, l <p\\j+k+l=p \endSb \frac {p!}{j!k!l!} x^jy^kz^l,$$ 
it is easily seen that all the coefficients of $H_p$ are divisible by $p$, and that $H_p(x,-x,z)=H_p(x,y,-y)=H_p(x,y,-x)=0$, so that $H_p$ is divisible by $p$, $x+y$, $y+z$ and $z+x$, in $Z[x,y,z]$, which are pairwise coprime. Therefore $H_p$ is divisible by their product $p(x+y)(y+z)(z+x)$, i.e., $E_3$ has its coefficients in $\Bbb Z$. Moreover, it is noted that 
$$E_3 = E_2 +Q(x,y,z),$$
where $Q \in z\Bbb Z[x,y,z]$. Indeed,
$$\align 
& E_3 - E_2 =  \frac {(x+y+z)^p -x^p -y^p -z^p}{p(x+y)(y+z)(z+x)} - \frac {(x+y)^p -x^p -y^p}{pxy(x+y)} = \\ 
&= \frac {xy((x+y+z)^p -x^p -y^p -z^p) -(y+z)(z+x)((x+y)^p -x^p -y^p)}{pxy(x+y)(y+z)(z+x)} = \\
&= \frac {xy \left(( \sum_{k=1}^{p-1} \binom pk (x+y)^k z^{p-k} +(x+y)^p -x^p -y^p \right) - (xy+xz+yz+z^2) ((x+y)^p -x^p-y^p)} {pxy(x+y)(y+z)(z+x)} = \\ 
& = \frac {xy  \sum_{k=1}^{p-1} \binom pk (x+y)^k z^{p-k} -z(x+y+z)  ((x+y)^p -x^p-y^p)}{pxy(x+y)(y+z)(z+x)} , 
\endalign$$
which is a polynomial in $\Bbb Z[x,y,z]$ all of whose terms are multiples of $z$. 

Then, came the main question. Setting $b=z(x+y+z)$, my correspondent made the observation that
$$E_3 = a_n +b(a_{n-1}+b(a_{n-2} + \cdots +b(a_2+b(a_1+b))\cdots )) = b^n +a_1b^{n-1}+\cdots a_{n-2}b^2 +a_{n-1}b +a_n,$$
where $n=\frac {p-3}2$ and all the $a_k \in \Bbb Z[x,y]$, with $a_n = E_2$. This amounts to stating that $E_3$ is a polynomial in $b$ with coefficients in $\Bbb Z[x,y]$, and asking for a proof. 

In April 2004, a second message came, extending the first one. In it, he makes the observation that $a_n +xy a_{n-1} = (x+y)^{p-3}$, and he conjectures that, for $0\leq m \leq n-1$, where $n=\frac {p-3}2$, 
$$a_{n-m}+ xya_{n-m-1} = a_{n-m}^{(1)} (x+y)^{p-2m-3},$$
where $a_k^{(j)}$ is the coefficient of the $j$-th term in the polynomial $a_k$.  Then it is claimed that a result of E. Catalan, cited in L. E. Dickson's \it {History of the Theory of Numbers, vol. II}\rm , gives the coefficients of $C_p(x,y)= \sum_{k=0}^{p-3 -2e} c_k x^{p-3-2e-k} y^k$, namely (with a correction)
$$c_k = \frac 1p \left(\binom {p-1}k -(-1)^k \right).$$
He then claims to deduce that $a_{n-1}^{(1)} = a_n^{(2)}$ and $a_n^{(2)} = \frac {p-3}2$. All computations were done, using \it {Mathematica}\rm. 

In May 2004, a third message came saying that he sent his query to Andrew Granville, who replied that any polynomial $F(x,y,z)$ which is invariant under the map 
$\varphi: z\mapsto -(x+y+z)$ (as is $E_3$) can be written as a polynomial in $b= z(x+y+z)$, providing a simple proof by induction on the degree of $F$. Indeed, setting $a_0 = F(x,y,0)$, Granville notes that $a_0$ and $F(x,y,z) -a_0$ are invariant under $\varphi$. But $F(x,y,z) -a_0$ is divisible by $z$, and, applying $\varphi$, it is also divisible by $-(x+y+z)$, and therefore it is divisible by $b$. Noting that $b$ is also invariant under $\varphi$, he deduces that $G(x,y,z) = \frac {F(x,y,z)-a_0}b$ is invariant under $\phi$, and  has lower degree than $F$. Thus, by the induction hypothesis, $G$ can be written as a polynomial in $b$, and therefore so can $F$.

It is this brilliantly simple argument of Granville that we attempt to generalize below.  

\medpagebreak
$\text{\bf Remarks.\ }$\rm The polynomials denoted by $E_2$ and $E_3$ will henceforth be designated by $E_{2, p}$ and $E_{3, p}$ to emphasize their dependence on $p$. 
Note that 
$$\align \frac {(x+y)^p - (x^p+y^p)}{x+y} &=(x+y)^{p-1} - \frac {x^p+y^p}{x+y} = \sum_{k=0}^{p-1} \binom {p-1}k x^k y^{p-1-k} -  \sum_{k=0}^{p-1} (-1)^k x^k y^{p-1-k} \\ &=  \sum_{k=0}^{p-1} \left( \binom {p-1}k -(-1)^k \right) x^k y^{p-1-k} = \sum_{k=1}^{p-2} \left( \binom {p-1}k -(-1)^k \right) x^k y^{p-1-k} , \endalign$$
since $\binom {p-1}k -(-1)^k =0$ for $k=0$ and $k=p-1$. Hence 
$$\align E_{2,p}(x,y) &= \frac {(x+y)^p - (x^p+y^p)}{pxy(x+y)} = \frac 1{pxy} \sum_{k=1}^{p-2} \left( \binom {p-1}k -(-1)^k \right) x^k y^{p-1-k} = \\
&= \sum_{k=1}^{p-2} \frac {\binom {p-1}k -(-1)^k}{p}  x^{k-1} y^{p-k-2}, \endalign$$
where the coefficients $\frac {\binom {p-1}k -(-1)^k}{p}$ are integers, for $1\leq k \leq p-2$, since 
$$\binom {p-1}k = \frac {(p-1)(p-2) \cdots (p-k)}{k!} \equiv \frac {(-1)^k k!}{k!} \equiv (-1)^k \pmod p.$$

On the other hand,  
$$\align H_p(x,y,z) &= (x+y+z)^p -(x^p+y^p+z^p) =\sum_{m=0}^p \binom pm (x+y)^m z^{p-m} - (x^p+y^p+z^p) = \\ &=\sum_{m=1}^p \binom pm (x+y)^m z^{p-m} - (x^p+y^p),  \endalign$$ 
so that 
$$\align E_{3,p}(x,y,z) &= \frac {H_p(x,y,z)}{p(x+y)(y+z)(z+x)} = \frac 1{p(y+z)(z+x)} \left( \sum_{m=1}^p \binom pm (x+y)^{m-1} z^{p-m} - \frac {x^p+y^p}{x+y} \right) \\ & = \frac 1{p(y+z)(z+x)} \left( \sum_{m=1}^p \binom pm (x+y)^{m-1} z^{p-m} -  \sum_{k=0}^{p-1} (-1)^k x^k y^{p-1-k} \right) = \\ & 
=  \frac 1{p(y+z)(z+x)} \sum_{k=0}^{p-1} \left( \binom p{k+1} (x+y)^k z^{p-k-1} - (-1)^k x^k y^{p-k-1} \right), \endalign $$
the last equality being obtained by setting $k=m-1$ in the first of the two sums. 

Furthermore, multiplying and dividing by $(x+y)$, distributing in the numerator, and setting $l=k+1$, we get 
$$\align E_{3,p}(x,y,z) &= \frac 1{p(x+y)(y+z)(z+x)} \left( \sum_{l=1}^p \binom pl (x+y)^l z^{p-l} + \sum_{l=1}^p (-1)^l x^l y^{p-l} - \sum_{k=0}^{p-1} (-1)^k x^k y^{p-k} \right)  \\ & = \frac 1{p(x+y)(y+z)(z+x)} \left( \sum_{l=0}^p \binom pl (x+y)^l z^{p-l} -x^p -y^p -z^p) \right). \endalign$$
Expanding the binomial $(x+y)^l$, and replacing all binomial coefficients by their defining expressions, then reversing the order of summations, and factoring out the common factors, we get 
$$E_{3,p}(x,y,z) = \frac 1{p(x+y)(y+z)(z+x)} \left( p! z^p\sum_{k=0}^p \frac 1{k!} \left( \frac xy \right)^k \sum_{l=k}^p\frac 1{(p-l)!(l-k)!} \left( \frac yz \right)^l -x^p -y^p -z^p \right) .$$
Hence the following two expressions for $E_{3,p}$, established above, 
$$\align E_{3,p}(x,y,z) &= \frac 1{p(y+z)(z+x)}  \sum_{k=0}^{p-1} \left( \binom p{k+1} (x+y)^k z^{p-k-1} - (-1)^k x^k y^{p-k-1} \right) = \\
&=  \frac 1{p(x+y)(y+z)(z+x)} \left( p! z^p\sum_{k=0}^p \frac 1 {k!} \left( \frac xy \right)^k \sum_{l=k}^p\frac 1{(p-l)!(l-k)!} \left( \frac yz \right)^l -x^p -y^p -z^p \right). \endalign$$

Similarly, exchanging $x+y$ and $y+z$, then $x+y$ and $x+z$, gives, respectively,
$$\align E_{3,p}(x,y,z) &= \frac 1{p(x+y)(x+z)} \sum_{k=0}^{p-1} \left( \binom p{k+1} (y+z)^k x^{p-k-1} - (-1)^k y^k z^{p-k-1} \right) = \\
&=  \frac 1{p(x+y)(y+z)(z+x)} \left( p! x^p\sum_{k=0}^p \frac 1 {k!} \left( \frac yz \right)^k \sum_{l=k}^p\frac 1{(p-l)!(l-k)!} \left( \frac zx \right)^l -x^p -y^p -z^p \right), \endalign$$
and
$$\align E_{3,p}(x,y,z) &= \frac 1{p(x+y)(y+z)} \sum_{k=0}^{p-1} \left( \binom p{k+1} (x+z)^k y^{p-k-1} - (-1)^k x^k z^{p-k-1} \right) = \\
&=  \frac 1{p(x+y)(y+z)(z+x)} \left( p! y^p\sum_{k=0}^p \frac 1 {k!} \left( \frac zx \right)^k \sum_{l=k}^p\frac 1{(p-l)!(l-k)!} \left( \frac xy \right)^l -x^p -y^p -z^p \right). \endalign$$

Moreover, after a search for the result of E. Catalan mentioned above, we located it in an article (\cite {1}), originally published in 1861, then taken up again, and complemented, in 1885, by Catalan, where he gives an explicit expression for the polynomials 
$$E_{3,n}(x,y,z) = \frac {(x+y+z)^n - x^n -y^n -z^n}{(x+y)(y+z)(z+x)},$$ 
for odd integers $n>3$. He defines the homogeneous polynomials of degree $m\geq 0$, with all coefficients equal to 1 
$$H_m = H_m(x,y,z) = \sum_{j,k,l \in \Bbb N: j+k+l = m} x^j y^k z^l = \sum_{j=0}^m \sum_{k=0}^{m-j} x^j y^k z^{m-j-k},$$
with $H_0 =1$, he sets $P = x+y+z$, and he establishes that, for odd $n>3$,
$$E_{3,n}(x,y,z) = \sum_{m=0}^{n-3} H_m(x,y,z) P^{n-3-m} + 2 H_{\frac {n-3}2} (x^2, y^2, z^2).$$

 \pagebreak

\head{The variation}  \endhead

\bigskip

Let $K$ be a commutative field. 
Take $n+1$ variables $x = (x_1,...,x_n)$ and $z$, then consider the  $K[x]$-algebra of  polynomials $A = K[x,z]$. Choose a polynomial, $p(x,z)$, in $A$ and define the following  transformation $T$ of $A$ into $A$: To each  polynomial $F(x,z)$, it matches the polynomial
$$TF(x,z) = F(x,p(x,z)).\tag{1}$$
This transformation $T$  is a $K[x]$-homomorphism of the algebra  $A$ into itself. Le $B$ be the set of  polynomial  invariant under $T$, i.e.,
$$B = \{F\in A : TF = F\}.\tag{2}$$
Of course, $B$ is a $K[x]$-subalgebra of $A$ and it contains $K[x]$. Then set
$$C = B \setminus K[x].\tag{3}$$
How can the elements of $B$ be characterized?

\noindent What is needed for $C$ not to be empty?

\noindent The case $p(x,z) = z$ is trivial since $T =$ \rm{id}$_A$ and $B = A$.

\

\noindent We  shall say that polynomial $p(x,z)$ is {\it \`a la Granville} when
$$p(x,z) = -z + r(x).\tag{4}$$
Notice that  $T$ is an involution (i.e., $T^2  =$ \rm{id}$_A$) in that case since we  will have $Tp(x,z) = z$.

\noindent The case studied by Granville is that in  which $p(x,z) = -(x+y+z)$.

\

\noindent Now, on to the  ``variation''!

\

\subhead{1 Preparation}\endsubhead  Introduce the sequence $T^k$ ot iterates of the transformation $T$, and the related sequence of polynomials $p_k(x,z) = T^kz$. 

\

\noindent Of course, $T^0$ is \rm{id}$_A$,  and $T^1 = T$.

\

\noindent Similarily, $p_0(x,z) = z$ and $p_1(x,z)$ is $p(x,z)$ itself.

\

\noindent For short, we say that $p(x,z)$ is   $1$-adequate when $p(x,z) = z$.

\

\noindent Let $m > 1$ be an integer. Henceforth, we assume that $K$ contains the $m$-th roots of unity. We say that a polynomial $p(x,z$) is $m$-adequate when it has the following form:
 $$p(x,z) = q z + r(x), \tag{*}$$ 
where $q$ is a primitive $m$-th root of unity, and $r(x)$ is in  $K[x]$.     

\

\noindent So, $p(x,z)$ is ``\`a la Granville'' if and only if it is  $2$-adequate.

\

\noindent When $p(x,z)$ is $m$-adequate for some $m \geq 1$, just say it is adequate.

\

\subhead{2  Remark}\endsubhead When $p(x,z)$ is $m$-adequate, we have $p_m(x,z) = z$ and  therefore  $T^m =$ \rm{id}$_A$. More precisely, we have

\

$p_0(x,z) = z$

$p_1(x,z) = q z + r(x)$

$p_2(x,z) = q^2 z+ (q + 1).r(x)$

$\ldots \ldots \ldots \ldots \ldots \ldots \ldots \ldots \ldots$

$p_k(x,z) = q^k z + (q^{k-1} + \dots + q + 1).r(x)$

$\ldots \ldots \ldots \ldots \ldots \ldots \ldots \ldots \ldots \ldots \ldots \ldots \ldots \ldots \ldots$

$p_m(x,z) = z$.

\

\noindent When $r(x) \neq 0$, polynomials $p_0, p_1,\dots,p_{m-1}$ are two by two mutually prime.

\

\noindent Let  $d(P)$ be the degree in  $z$ of each polynomial $P(x,z)$.]

\

\subhead{3 Theorem}\endsubhead Let $p(x,z)$ be $m$-adequat. Let $b$ be the product of polynomials $\{ p_k  :  0 \leq k <m \}$. A polynomial $F(x,z)$ is in $B$ if and only if it has the following form:
 $$F(x,z) = G(b), \ \text{with} \ G(y)\in  K[x][y].\tag{6}$$

\

\subhead{Proof}\endsubhead Observe that polynomial $b$ is  in $B$, so $TG(b) = G(b)$, and every polynomial which has the form (6) is then in $B$.

\

\noindent For the converse, distinguish two cases.

\

\noindent If $r(x) = 0$, then $F$ is in $B$ means $F(x,z) = F(x, q z)$; since $q^m = 1$, we conclude that $F$ is a  polynomial in $z^m$.

\

\noindent Otherwise, proceed by induction on the degree, $d(F)$, of $F$ in $z$.

\

\noindent The case $d(F) = 0$ is trivial.    

\

\subhead{In four steps}\endsubhead

\

\subhead{Lemma 1}\endsubhead If $F(x,z)$ is in $B$, each of the $p_k(x,z)$, $k \geq 0$, divides the polynomial $W(x,z):= F(x,z) - F(x,0)$.

\

\noindent Indeed, $W(x,z)$ also is in $B$ and  is divisible by $z$. Therefore, we must have $W(x,z) = z.U(x,z)$, and for $k\geq 1$, we have 

$W(x,z) = T^kU(x,z) = p_k(x,z).T^kU(x,z)$.

\noindent If $F(x,z)$ is in $B$, polynomial $F(x,z) - F(x,0)$ is  divisible  by $b$, i.e., we have $F(x,z) - F(x,0) = b(x,z).V(x,z)$. So, $V(x,z)$ is in $B$, moreover  we have $d(V) < d(F)$. Whence the conclusion.\qed

\

\subhead{Complement}\endsubhead The following three conditions are equivalent.

\

\ \ (i) The set $C$ is not empty.

\ (ii)  There exists an integer  $m \geq 1$ such that $T^m =$ \rm{id}$_A$.

(iii) Polynomial $p(x,z)$ is adequate.

\

\subhead{Proof}\endsubhead If $p(x,z)$ is $m$-adequate, we already know that $C$ is not empty, and that $T^m =$ \rm{id}$_A$.\qed

\

\subhead{Lemma  2}\endsubhead If $d(p) \neq 1$, then $B = K[x]$.

\

\noindent Indeed, $d(TF) = d(F).d(p)$. If  $F$ is in $B$, i.e., $TF = F$, we have $d(F) = d(F).d(p)$, and if, moreover, $d(p)\neq 1$,  then $d(F) = 0$. \qed

\

\subhead{Lemma 3}\endsubhead If $C$ is not empty, the polynomial $p(x,z)$ is adequate.

\subhead{Proof}\endsubhead By the previous lemma, we know that $p(x,z)$ has the following form:
 $$p(x,z) = q(x).z + r(x), \ \text{with} \ q(x)\neq 0. \tag{7}$$
Take $F \in B$ such that $d(F) = m \geq 1$.

\noindent Ordering $F$ in decreasing powers of  $z$, write

$F = u.z^m + v.z^(m-1) + \cdots$

\noindent whence

$TF = u.(q z + r)^m + v.(q z + r)^{m-1} + \cdots$

    $= u.q^m z^m + (u.m.q^{m-1}.r + v.q^{m-1}z^{m-1} + \cdots$

\noindent then,  identifying the first two coefficients, obtain both equalities

 $$q^m = 1\tag{8}$$

$$q^{m-1}.(m.r + v) = v.\tag{9}$$

\noindent For $q = 1$, we therefore have  $r = 0$. \qed

\subhead{Lemma 4}\endsubhead If  $T^m =$ \rm{id} $_A$  for some $m \geq 1$, the product $b$ of  polynomials

\noindent  in $\{ p_k  :  0 \leq k < m\}$ is invariant by  $T$, and $d(b) \geq 1$.

\

\noindent {\bf This ends the proof}\hfill{\;\;\;\qed\qed}

\

\enddocument